\documentclass[pre,twocolumn,showpacs,superscriptaddress,preprintnumbers,floatfix]{revtex4}

\usepackage{amssymb,amsmath} 
\usepackage{graphicx}
\usepackage{bm}
\usepackage{vaucanson-g} 

\newcommand{\eqnref}[1]{Eq.~(\ref{#1})}
\newcommand{\begeqnref}[1]{Equation~\ref{#1}}
\newcommand{\figref}[1]{Fig.~\ref{#1}}
\newcommand{\begfigref}[1]{Figure~\ref{#1}}
\newcommand{\appref}[1]{App.~\ref{#1}}

\newcommand{\hmu}	{h_\mu} 
\newcommand{\KLd}	{\mathcal{D}} 
\newcommand{\EQP}	{E(Q,P)}
\newcommand{\Asize}[0]{\vert \mathcal{A} \vert}
\newcommand{\hk}[0]{ \overleftarrow{s}^k} 
\newcommand{\htl}[2]{ \overleftarrow{s}_{#1}^{#2}} 

\newcommand{\nsk}[0]{n(\overleftarrow{s}^k)} 
\newcommand{\nsks}[0]{n(\overleftarrow{s}^k s)} 
\newcommand{\ask}[0]{\alpha(\overleftarrow{s}^k)} 
\newcommand{\asks}[0]{\alpha(\overleftarrow{s}^k s)} 
\newcommand{\MP}[0]{\mathbf{\theta}} 
\newcommand{\MPk}[0]{\mathbf{\theta}_k} 
\newcommand{\MC}[0]{\mathbf{M}} 
\newcommand{\MCk}[0]{\mathbf{M}_k} 
\newcommand{\MCkprime}[0]{\mathbf{M}_{k}'} 

\newcommand{\avg}[2]{\mathbf{E}_{#2}[#1]} 
\newcommand{\var}[2]{\mathbf{Var}_{#2}[#1]} 
\newcommand{\cov}[2]{\mathbf{Cov}[#1,#2]} 
\newcommand{\psk}[0]{p(\overleftarrow{s}^k)}  
\newcommand{\psks}[0]{p(s\vert \overleftarrow{s}^k)}  
\newcommand{\qsk}[0]{q(\overleftarrow{s}^k)}  
\newcommand{\qsks}[0]{q(s\vert \overleftarrow{s}^k)}
\newcommand{\rsk}[0]{r(\overleftarrow{s}^k)}  
\newcommand{\rsks}[0]{r(s\vert \overleftarrow{s}^k)}  
\newcommand{\usks}[0]{u(s\vert \overleftarrow{s}^k)}  
\newcommand{\usk}[0]{u(\overleftarrow{s}^k)}  

\newcommand{\mle}[0]{MLE} 
\newcommand{\pme}[0]{PME} 

\begin{document}

\preprint{Santa Fe Institute Working Paper 07-03-XXX}
\preprint{arxiv.org/xxxxx/0703XXX}

\title{Inferring Markov Chains: Bayesian Estimation,\\
Model Comparison, Entropy Rate, and Out-of-class Modeling}

\author{Christopher~C.~Strelioff}
	\email{streliof@uiuc.edu} 
	\affiliation{Center for Computational Science \& 
	Engineering and Physics Department,\\ 
	University of California at Davis, One Shields Avenue, Davis, CA 95616}
 	\affiliation{Center for Complex Systems Research 
 	and Physics Department,\\
	University of Illinois at Urbana-Champaign, 
	1110 West Green Street, Urbana, Illinois 61801}
\author{James P. Crutchfield}
	\email{chaos@cse.ucdavis.edu}
	\affiliation{Center for Computational Science \& 
	Engineering and Physics Department,\\ 
	University of California at Davis, One Shields Avenue, Davis, CA 95616}
\author{Alfred W. H\"{u}bler}
	\email{a-hubler@uiuc.edu}
 	\affiliation{Center for Complex Systems Research 
 	and Physics Department,\\
	University of Illinois at Urbana-Champaign, 
	1110 West Green Street, Urbana, Illinois 61801}

\begin{abstract} 
Markov chains are a natural and well understood tool for describing
one-dimensional patterns in time or space. We show how to infer $k$-th order
Markov chains, for arbitrary $k$, from finite data by applying Bayesian
methods to both parameter estimation and model-order selection. Extending
existing results for multinomial models of discrete data, we connect inference
to statistical mechanics through information-theoretic (type theory) techniques.
We establish a direct relationship between Bayesian evidence and the partition
function which allows for straightforward calculation of the expectation and
variance of the conditional relative entropy and the source entropy rate.
Finally, we introduce a novel method that uses finite data-size scaling with
model-order comparison to infer the structure of out-of-class processes.
\end{abstract}

\pacs{02.50.Tt,02.50.Ga,05.10.Gg}
                         
\maketitle

%
%
\section{Introduction}

Statistical inference of models from small data samples is a vital tool in
the understanding of natural systems.  In many problems of interest data
consists of a sequence of \emph{letters} from a finite \emph{alphabet}.  
Examples include analysis of sequence information in
biopolymers~\cite{Avery1999,JSLiu1999}, investigation of
one-dimensional spin systems~\cite{Crutchfield1997}, models of natural 
languages~\cite{MacKay1994}, and coarse-grained models of chaotic
dynamics~\cite{Crutchfield1983,BLHao1998}.  This diversity of potential
application has resulted in the development of a variety of representations
for describing discrete-valued data series.

We consider the $k$-th order Markov chain model class which uses the previous
$k$ letters in a sequence to predict the next letter. Inference of Markov
chains from data has a long history in mathematical statistics.  Early work
focused on maximum likelihood methods for estimating the parameters of the
Markov chain~\cite{TWAnderson1957,Billingsley1961a,Chatfield1973}. This work 
often assumed a given fixed model order. That is, no \emph{model comparison}
across orders is done. This work also typically relied on the assumed
asymptotic normality of the likelihood when estimating regions of
confidence and when implementing model comparison.  As a result, the realm
of application has been limited to data sources where these conditions are
met.  One consequence of these assumptions has been that data sources which
exhibit \emph{forbidden words}, symbol sequences which are not allowed, cannot
be analyzed with these methods.  This type of data violates the assumed
normality of the likelihood function. 

More recently, model comparison in the maximum likelihood approach has been
extended using various \emph{information criteria}. These methods for
model-order selection are based on extensions of the likelihood ratio and allow
the comparison of more than two candidate models at a time. The most widely used
are \emph{Akaike's information criteria} (AIC)~\cite{HTong1975} and the
\emph{Bayesian information criteria} (BIC)~\cite{Katz1981}. (Although the
latter is called Bayesian, it does not employ Bayesian model comparison in
the ways we will present here.) In addition to model selection using information
criteria, methods from information theory and machine learning have also been
developed.  Two of the most widely employed are \emph{minimum
description length} (MDL)~\cite{JRissanen1984} and \emph{structural risk
minimization}~\cite{VVapnik1999}.  Note that MDL and Bayesian
methods obtain similar results in some situations~\cite{Vitanyi2000}.  However,
to the best of our knowledge, structural risk minimization has not been adapted
to Markov chain inference.

We consider Bayesian inference of the Markov chain model class, extending
previous results~\cite{MacKay1994,JSLiu1999,Baldi2001,Durbin1998}. We provide
the details necessary to infer a Markov chain of arbitrary order, choose
the appropriate order (or weight orders according to their probability),
and estimate the data source's entropy rate.  The latter is important for
estimating the intrinsic randomness and achievable compression rates for
an information source~\cite{Cover1991}.  The ability to weight Markov chain
orders according their probability is unique to Bayesian methods and
unavailable in the model selection techniques discussed above.

In much of the literature just cited, steps of the inference process
are divided into (i) point estimation of model parameters, (ii) model
comparison (hypothesis testing), and (iii) estimation of functions of the
model parameters. Here we will show that Bayesian inference connects all
of these steps, using a unified set of ideas. Parameter estimation is the first
step of inference, model comparison a second level, and estimation of the
entropy rate a final step, intimately related to the mathematical structure
underlying the inference process.  This view of connecting model to data
provides a powerful and unique understanding of inference not available in the
classical statistics approach to these problems. As we demonstrate, each of
these steps is vital and implementation of one step without the others does
not provide a complete analysis of the data-model connection. 

Moreover, the combination of inference of model parameters, comparison of
performance across model orders, and estimation of entropy rates provides a
powerful tool for understanding Markov chain models themselves. Remarkably,
this is true even when the generating data source is outside of the Markov
chain model class.
Model comparison provides a sense of the structure of the data source, whereas
estimates of the entropy rate provide a description of the inherent randomness. 
Bayesian inference, information theory, and tools from statistical mechanics
presented here touch on all of these issues within a unified framework.

We develop this as follows, assuming a passing familiarity with Bayesian
methods and statistical mechanics. First, we discuss estimation of Markov
chain parameters using Bayesian methods, emphasizing the use of the complete
marginal posterior density for each parameter, rather than point estimates
with error bars. Second, we consider selection of the appropriate memory
$k$ given a particular data set, demonstrating that a mixture of orders may
often be more appropriate than selecting a single order. This is certainly
a more genuinely Bayesian approach. In these first two parts
we exploit different forms of Bayes' theorem to connect data and model class.

Third, we consider the mathematical structure of the evidence (or marginal
likelihood) and draw connections to statistical mechanics.  In this discussion
we present a method for estimating entropy rates by taking derivatives of a
partition function formed from elements of each step of the inference procedure.
Last, we apply these tools to three example information sources of increasing
complexity. The first example belongs to the Markov chain model class, but
the other two are examples of hidden Markov models (HMMs) that fall outside
of that class. We show that the methods developed here provide a powerful tool
for understanding data from these sources, even when they do not belong to the
model class being assumed.

\section{Inferring Model Parameters}

In the first level of Bayesian inference we develop a systematic relation
between the data $D$, the chosen \emph{model class} $M$, and the vector of
\emph{model parameters} $\MP$. The object of interest in the inference of
model parameters is the \emph{posterior probability density}
$P\left( \MP \vert D, M \right)$.  This is the probability of the model
parameters given the observed data and chosen model. To find the posterior
we first consider the joint distribution $P\left( \MP, D \vert M \right)$
over the data and model parameters given that one has chosen to model the
source with a representation in a certain class $M$. This can be factored in
two ways: $P\left( \MP \vert D, M \right)P\left(D \vert M\right)$ or
$P\left( D \vert \MP, M \right)P\left(\MP \vert M\right)$.  Setting these
equal and solving for the posterior we obtain Bayes' theorem:
\begin{equation}
\label{eqn:bayes}
P\left( \MP \vert D, M \right) 
	= \frac{ P\left( D \vert \MP , M \right) \; 
	P\left( \MP \vert M \right) }{ P\left( D \vert M \right) }.
\end{equation}

The \emph{prior} $P\left( \MP \vert M \right)$ specifies our assumptions
regarding the model parameters. We take a pragmatic view of the prior,
considering its specification to be a statement of assumptions about the
chosen model class. The \emph{likelihood} $P\left( D \vert \MP , M \right)$
describes the probability of the data given the model.  Finally, the
\emph{evidence} (or marginal likelihood) $P\left( D \vert M \right)$ is the
probability of the data given the model.  In the following sections we
describe each of the quantities in detail on our path to giving an explicit
expression for the posterior.

\subsection{Markov chains}

The first step in inference is to clearly state the assumptions that make up
the model.  This is the foundation for writing down the likelihood of a data
sample and informs the choice of prior. We assume that a single data set of
length $N$ is the starting point of the inference and that it consists of
\textit{symbols} $s_t$ from a finite alphabet $\mathcal{A}$,  
\begin{equation}
	\label{eqn:data}
	D = s_0 s_1 \ldots s_{N-1} \; , \; s_t \in \mathcal{A}.
\end{equation}
We introduce the notation $\htl{t}{k}$ to indicate a length-$k$ sequence of
letters ending at position $t$: e.g., $\htl{4}{2}=s_3s_4$.

The $k$-th order Markov chain model class assumes finite memory and 
stationarity in the data source.  The finite memory condition, a
generalization of the conventional Markov property, can be written
\begin{equation}
p(D)	 = p(\htl{k-1}{k}) \prod_{t=k-1}^{N-2} p(s_{t+1} \vert \htl{t}{k}) ~,
	    \label{eqn:markov_condition}
\end{equation}
thereby factoring into terms which depend only on preceding words of
length-$k$. The stationarity condition can be expressed
\begin{equation}
	\label{eqn:stationarity}
	p(s_t \vert \htl{t-1}{k}) = p(s_{t+m} \vert \htl{t+m-1}{k}) ~, 
\end{equation}
for any $(t,m)$.  \begeqnref{eqn:stationarity} results in a simplification of 
the notation because we no longer need to track the position index, 
$p(s_t = s \vert \htl{t-1}{k} = \hk ) = p( s \vert \hk )$ for any $t$.  Given 
these two assumptions, the model parameters of the $k$-th order Markov chain
$\MCk$ are
\begin{equation}
	\label{eqn:model_parameters}
 	\MPk  = \left\{ \, p( s \vert \hk ) : s \in \mathcal{A}, 
 	\hk \in \mathcal{A}^k \, \right\}.
\end{equation}
A normalization constraint is placed on these parameters $\sum_{s\in
\mathcal{A}} p( s \vert \hk ) = 1$ for each word $\hk$.

The next step is to write down the elements of Bayes' theorem specific to the
$k$-th order Markov chain.

\subsection{Likelihood}

Given a sample of data $D=s_{0}s_{1} \ldots s_{N-1}$, the likelihood can be 
written down using the Markov property of~\eqnref{eqn:markov_condition} and the 
stationarity of~\eqnref{eqn:stationarity}.  This results in the form
\begin{equation}
	\label{eqn:likelihood}
	P(D\vert \MPk, \MCk) = \prod_{ s \in \mathcal{A} } 
	\prod_{ \hk \in \mathcal{A}^{k} } p( s \vert \hk )^{\nsks} ,
\end{equation}
where $\nsks$ is the number of times the \textit{word} $\hk s$ occurs in the
sample $D$.  For future use we also introduce notation for the number of times a
word $\hk$ has been observed $\nsk = \sum_{s \in \mathcal{A}} \nsks$.  We note
that~\eqnref{eqn:likelihood} is conditioned on the \emph{start sequence}
$\hk = s_0s_1\ldots s_{k-1}$.

\vspace{-0.125in}
\subsection{Prior}
\vspace{-0.125in}

The prior $P(\theta|M)$ is used to specify assumptions about the model to be
inferred before the data is considered. Here we use
\emph{conjugate priors} for which the posterior distribution has the same
functional form as the prior.  Our choice allows us to derive exact expressions
for many quantities of interest in inference. This provides a powerful tool for
understanding what information is gained during inference and,
especially, model comparison.

The exact form of the prior is determined by our assignment of 
\emph{hyperparameters} $\asks$ for the prior which balance the strength of
the modeling assumptions encoded in the prior against the weight of the data.
For a $k$-th order Markov chain, there is one hyperparameter for each word
$\hk s$, given the alphabet under consideration. A useful way to think about
the assignment of values to the
hyperparameters is to relate them to fake counts $\tilde{n}(\hk s)$, such that
$\asks = \tilde{n}(\hk s) + 1$.  In this way, the $\asks$ can be set to reflect
knowledge of the data source and the strength of these prior assumptions can be
properly weighted in relation to the actual data counts $\nsks$.

The conjugate prior for Markov chain inference is a product of Dirichlet
distributions, one for each word $\hk$. It restates the finite-memory
assumption from the model definition:
\begin{eqnarray}
	P(\MPk \vert \MCk ) 
	& = & \prod_{\hk \in \mathcal{A}^{k}} \left\{
	\frac{ \Gamma( \ask  )}{
	\prod_{s\in\mathcal{A}} \Gamma( \asks ) } \right. \nonumber \\
    & \times & \delta \mathbf{(}1-\sum_{s\in\mathcal{A}} 
    p( s \vert \hk )\mathbf{)} \label{eqn:prior} \\
	& \times & \left. \prod_{s\in\mathcal{A}} p( s \vert \hk )^{\asks-1} 
	\right\}. \nonumber
\end{eqnarray}
(See App. \ref{app:Dirichlet} for relevant properties of Dirichlet
distributions.)
The prior's hyperparameters $\{ \asks \}$ must be real and positive.  We
also introduce the more compact notation $\ask = \sum_{s \in \mathcal{A}}
\asks$.  The function $\Gamma(x)=(x-1)!$ is the well known Gamma function.  The
$\delta$-function constrains the model parameters to be properly normalized:
$\sum_{s \in \mathcal{A}} \psks = 1$ for each $\hk$.

Given this functional form, there are at least two ways to interpret what the
prior says about the Markov chain parameters $\MPk$. In addition to considering
fake counts $\tilde{n}( \cdot )$, as discussed above, we can consider the
range of fluctuations in the estimated $\psks$. Classical statistics would
dictate describing the fluctuations via a single value with error bars. This
can be accomplished by finding the average and variance of $\psks$ with
respect to the prior. The result is:
\begin{eqnarray}
	\label{eqn:prior_mean}
	\avg{\psks}{\rm{prior}} & = & \frac{\asks}{\ask}~, \\
	\label{eqn:prior_variance}
	\var{\psks}{\rm{prior}} & = & \frac{\asks(\ask-\asks)}{\ask^2(1+\ask)} .
\end{eqnarray}

A second method, more in line with traditional Bayesian estimation, is to
consider the marginal distribution for each model parameter. For a Dirichlet
distribution, the marginal for any one parameter will be a Beta distribution.
With this knowledge, a probability density can be provided for each Markov chain
parameter given a particular setting for the hyperparameters $\asks$. In this
way, the prior can be assigned and analyzed in substantial detail.

A common stance in model inference is to assume all things are a-priori
equal.  This can be expressed by assigning $\asks=1$ for all $\hk \in
\mathcal{A}^k$ and $s \in \mathcal{A}$, adding \textit{no} fake counts
$\tilde{n}(\hk s)$.  This assignment results in a uniform prior distribution
over the model parameters and a prior expectation:
\begin{equation}
\avg{p(s\vert \hk)}{\rm{prior}} = 1/ \vert \mathcal{A} \vert ~.
\end{equation}

\vspace{-0.20in}
\subsection{Evidence}
\vspace{-0.125in}

Given the likelihood and prior derived above, the evidence $P(D|M)$ is seen
to be a simple normalization term in Bayes' theorem.  In fact, the evidence
provides the probability of the data given the model $\MCk$ and so plays a
fundamental role in model comparison.  Formally, the definition is
\begin{equation}
	P(D\vert \MCk ) = 	\int \; d\MPk \; P(D\vert \MPk, \MCk) 
						P(\MPk \vert \MCk ),
	\label{eqn:evidence_defn}
\end{equation}
where we can see that this term can be interpreted as an average of the 
likelihood over the prior distribution.  Applying this to the likelihood
in~\eqnref{eqn:likelihood} and the prior in~\eqnref{eqn:prior} produces 
\begin{eqnarray}
	P(D\vert \MCk) & = & \prod_{\hk \in \mathcal{A}^{k}} \left\{ \; 
	\frac{ \Gamma(\ask) }{ \prod_{s\in \mathcal{A}} \Gamma(\asks)} 
	\right. \nonumber \\
	& & \label{eqn:evidence} \\
	& \times & \left. 
	\frac{ \prod_{s\in \mathcal{A}} \Gamma(\nsks+\asks) }{ \Gamma(\nsk+\ask) }
	\; \right\}. \nonumber
\end{eqnarray}
As we will see, this analytic expression results in the ability to make useful
connections to statistical mechanics techniques when estimating entropy rates.
This is another benefit of choosing a conjugate prior with known properties.

\subsection{Posterior}

Using Bayes' theorem~\eqnref{eqn:bayes} the results of the three previous 
sections can be combined to obtain the posterior distribution over the
parameters of the $k$-th order Markov chain. One finds:
\begin{eqnarray}
	P(\MPk\vert D, \MCk) & = & \prod_{\hk \in \mathcal{A}^{k}} \left\{
	\frac{ \Gamma( \nsk + \ask  ) }{
	\prod_{s\in\mathcal{A}} \Gamma( \nsks + \asks ) } \right. \nonumber \\
	& \times & \delta \mathbf{(}1-\sum_{s\in\mathcal{A}} p( s \vert \hk)
	\mathbf{)} \label{eqn:posterior} \\
	& \times & \left. \prod_{s\in\mathcal{A}} 
	p( s \vert \hk )^{\nsks + \asks - 1} \right\}. \nonumber
\end{eqnarray}
As noted in selecting the prior, the resulting form is a Dirichlet 
distribution with modified parameters.  This is a result of choosing the 
conjugate prior: cf. the forms of \eqnref{eqn:prior} and
\eqnref{eqn:posterior}.

From~\eqnref{eqn:posterior} the estimation of the model parameters 
$p(s\vert \hk)$ and the uncertainty of these estimates can be given using the
known properties of the Dirichlet distribution.  As with the prior,
there are two main ways to understand what the posterior tells us about the
fluctuations in the estimated Markov chain parameters. The first uses a point
estimate with ``error bars''. We obtain these from the mean and variance of
the $\psks$ with respect to the posterior, finding
\begin{gather}
 	\avg{p(s\vert \hk)}{\rm{post}} =  \frac{ \nsks + \asks }{ \nsk + \ask }
	\label{eqn:posterior_mean} ~, \\ \nonumber \\
 	\var{p(s\vert \hk)}{\rm{post}}  =  \frac{ \nsks +\asks }{ ( \nsk + \ask )^2 }
	\nonumber  \\ \label{eqn:posterior_variance} \\
	\times  \frac{ ( \nsk + \ask ) - ( \nsks + \asks  ) }{
	( \nsk + \ask +1 ) }. \nonumber
\end{gather}
This is the \textit{posterior mean estimate} (\pme) of the model parameters.

A deeper understanding of~\eqnref{eqn:posterior_mean} is obtained through a
simple factoring:
\begin{eqnarray}
	\avg{p(s\vert \hk)}{\rm{post}} & = & \frac{1}{ \nsk + \ask } 
	\left[ \nsk \, \left (\frac{\nsks}{\nsk} \right) \right. \nonumber \\
	&& 	\label{eqn:pme_factor} \\
	& + & \left. \ask \, \left(\frac{\asks}{\ask} \right) \right], \nonumber
\end{eqnarray}
where $\nsks /\nsk $ is the \emph{maximum likelihood estimate} (\mle) 
of the model parameters and $\asks /\ask$ is the prior expectation given
in~\eqnref{eqn:prior_mean}.  In this form, it is
apparent that the posterior mean estimate is a weighted sum of the \mle~and 
prior expectation.  As a result, we can say that the posterior mean and
maximum likelihood estimates converge to the same value for $\nsk \gg \ask$. 
Only when the data is scarce, or the prior is set with strong conviction, 
does the Bayesian estimate add corrections to the \mle.

A second method for analyzing the resulting posterior density is to consider the
marginal density for each parameter.  As discussed with the prior, the marginal
for a Dirichlet is a Beta distribution.  As a result, we can either provide
regions of confidence for each parameter or simply inspect the density function.
The latter provides much more information about the inference being made than
the point estimation just given.  In our examples, to follow shortly, we
plot the marginal posterior density for various parameters of interest
to demonstrate the wealth of information this method provides.

Before we move on, we make a final point regarding the estimation of inference 
uncertainty. The form of the posterior is not meant to reflect the potential
fluctuations of the data source.  Instead, the width of the distribution
reflects the possible Markov chain parameters which are consistent with
observed data sample.  These are distinct notions and should not be conflated.

\subsection{Predictive distribution}

Once we have an inferred model, a common task is to estimate the probability of
a new observation $D^{(new)}$ given the previous data and estimated model.
This is implemented by taking an average of the likelihood of the new data:
\begin{equation}
P(D^{(new)}\vert \MPk, \MCk)
  = \prod_{\hk \in \mathcal{A}^k, s \in \mathcal{A}} p(s\vert \hk)^{m(\hk s)}
\end{equation}
with respect to the posterior
distribution~\cite{MacKay2003}:
\begin{eqnarray}
	\label{eqn:predictive_distribution_defn}
	P(D^{(new)}\vert D,\MCk) & =  & \int  d\MPk  P(D^{(new)}\vert \MPk, \MCk) \\ 
	& \times & P(\MPk \vert D, \MCk) ~. \nonumber
\end{eqnarray}
We introduce the notation $m(\hk s)$ to indicate the number of times the word
$\hk s$ occurs in $D^{(new)}$. This method has the desirable property, compared
to point estimates, that it takes into account the uncertainty in the model
parameters $\MPk$ as reflected in the form of the posterior distribution.

The evaluation of~\eqnref{eqn:predictive_distribution_defn} follows the same 
path as the calculation for the evidence and produces a similar
form; we find:
\begin{gather}
	P(D^{(new)}\vert D, \MCk)  =  \prod_{\hk \in \mathcal{A}^{k}} \left\{ \; 
	\frac{ \Gamma( \nsk+\ask) }{ \prod_{s\in \mathcal{A}} 
	\Gamma( \nsks + \asks)} \right. \nonumber \\
	\label{eqn:predictive_distribution} \\
	\times  \left. \frac{ \prod_{s\in \mathcal{A}} 
	\Gamma( \nsks + m(\hk s) + \asks ) }{ \Gamma( \nsk + m(\hk) + \ask ) }
	\; \right\}. \nonumber
\end{gather}

\section{Model Comparison}

With the ability to infer a Markov chain of a given order $k$, a common sense 
question is to ask how do we choose the correct order given a particular data 
set?  Bayesian methods have a systematic way to address this through 
the use of \emph{model comparison}.

In many ways, this process is analogous to inferring model parameters
themselves, which we just laid out.  We start by enumerating the set of model
orders to be compared $\mathcal{M} = \{ \MCk \}_{k_{min}}^{k_{max}}$, where
$k_{min}$ and $k_{max}$ correspond to the minimum and maximum order to be
inferred, respectively.  Although we will not consider an independent,
identically distributed (IID) model ($k=0$) here, we do note that this could
be included using the same techniques described below.

We start with the joint probability $P(M_{k},D \vert \mathcal{M} )$ of a
particular model $M_{k} \in \mathcal{M}$ and data sample $D$, factoring it in
two ways following Bayes' theorem. Solving for the probability of a particular
model class we obtain
\begin{equation}
	\label{eqn:model_comparison}
 	P(\MCk \vert D , \mathcal{M} ) = \frac{ P(D \vert \MCk, \mathcal{M} )
 					P(\MCk \vert \mathcal{M} ) }{ P(D \vert \mathcal{M})} ,	
\end{equation} 
where the denominator is the sum given by
\begin{equation}
P(D \vert \mathcal{M}) = 
  \sum_{\MCkprime \in \mathcal{M}}
  P(D \vert \MCkprime, \mathcal{M} )P(\MCkprime \vert \mathcal{M} ) ~.
\end{equation}
The probability of a particular model class in the set under consideration is
driven by two components: the evidence $P(D \vert \MCk, \mathcal{M})$, derived
in \eqnref{eqn:evidence}, and the prior over model classes
$P(\MCk \vert \mathcal{M} )$.

Two common priors in model comparison are: (i) all models are equally likely
and (ii) models should be penalized for the number of free parameters used to
fit the data.  In the first instance 
$P(\MCk \vert \mathcal{M})=1/ \vert \mathcal{M} \vert$ is the same for all 
orders $k$.  However, this factor cancels out because it appears in both the
numerator and denominator.  As a result, the probability of models using this
prior becomes
\begin{equation}
	\label{eqn:best_model_uniform_prior}
	P(\MCk \vert D , \mathcal{M} ) = \frac{P(D \vert \MCk, \mathcal{M} )
					}{
					\sum_{\MCkprime \in \mathcal{M}}
					P(D \vert \MCkprime, \mathcal{M} )}.
\end{equation}

In the second case, a common penalty for the number of model parameters is 
\begin{equation}
\label{eqn:df_penalty_prior}
P(\MCk \vert \mathcal{M}) = \frac{\exp( - \vert \MCk \vert )
						  }{\sum_{\MCkprime \in \mathcal{M}} 
						  \exp( - \vert \MCkprime \vert ) } ~,
\end{equation}
where $\vert \MCk \vert$ is the number of free parameters in the model. For a
$k$-th order Markov chain, the number of free parameters is
\begin{equation}
\vert \MCk \vert = \vert \mathcal{A} \vert^k(\vert \mathcal{A} \vert-1) ~,
\end{equation}
where $\vert \mathcal{A} \vert$ is the alphabet size. Thus, model
probabilities under this prior take on the form
\begin{equation}
	\label{eqn:best_model_df_penalty_prior}
	P(\MCk \vert D , \mathcal{M} ) = \frac{
					P(D \vert \MCk, \mathcal{M} ) \exp( - \vert \MCk \vert )
					}{
					\sum_{\MCkprime}
					P(D \vert \MCkprime, \mathcal{M} )
					\exp( - \vert \MCkprime \vert ) }.
\end{equation}
We note that the normalization sum in~\eqnref{eqn:df_penalty_prior}
cancels because it appears in both the numerator and denominator.

Bayesian model comparison has a natural \emph{Occam's razor} in the model
comparison process~\cite{MacKay2003}.  This means there is a natural preference
for smaller models even when a uniform prior over model orders is applied.  In
this light, a penalty for the number of model parameters can be seen as a very
cautious form of model comparison.  Both of these priors,
\eqnref{eqn:best_model_uniform_prior} and  
\eqnref{eqn:best_model_df_penalty_prior}, will be considered in
the examples to follow.

A note is in order on computational implementation. In general, the resulting
probabilities can be extremely small, easily resulting in numerical underflow
if the equations are not implemented with care. As mentioned
in~\cite{Durbin1998}, computation with extended logarithms can be used to
alleviate these concerns.

\section{Information Theory, Statistical Mechanics, and Entropy Rates}

An important property of an information source is its \emph{entropy rate}
$\hmu$, which indicates the degree of intrinsic randomness and controls the
achievable compression. A first attempt at estimating a source's entropy rate
might consist of plugging a Markov chain's estimated model parameters into the
known expression~\cite{Cover1991}. However, this does not
accurately reflect the posterior distribution derived above. This observation
leaves two realistic alternatives. The first option is to sample model
parameters from the posterior distribution. These samples can then be used to
calculate a set of entropy rate estimates that reflect the underlying posterior
distribution. A second option, which we take here, is to adapt methods from
type theory and
statistical mechanics previously developed for IID models~\cite{Samengo2002}
to Markov chains. To the best of our knowledge this is the first time these
ideas have been extended to inferring Markov chains; although cf.
\cite{Young1994}. 

In simple terms, type theory shows that the probability of an observed sequence
can be written in terms of the \emph{Kullback-Leibler} (KL) \emph{distance} and
the entropy rate.  When applied to the Markov chain inference problem the resulting
form suggests a connection to statistical mechanics. For example, we will show
that averages of the KL-distance and entropy rate with respect to the posterior
are found by taking simple derivatives of a partition function.  

The connection between inference and information theory starts by considering
the product of the prior~\eqnref{eqn:prior} and
likelihood~\eqnref{eqn:likelihood}:
\begin{equation}
P(\MPk\vert \MCk)P( D\vert \MPk, \MCk)=P( D, \MPk\vert \MCk) ~.
\end{equation}
This forms a joint distribution over the observed data $D$ and model parameters
$\MPk$ given the model order $\MCk$. Denoting the normalization constant from
the prior as $Z$ to save space, this joint distribution is
\begin{equation}
	\label{eqn:product_prior_likelihood}
 	P( D, \MPk\vert \MCk) = Z \, \prod_{\hk, s}  
 	p( s \vert \hk )^{\nsks + \asks - 1}.
\end{equation}
This form can be written, without approximation, in terms of conditional 
relative entropies $\KLd [\cdot \| \cdot ]$ and entropy rate $\hmu [\cdot]$:
\begin{eqnarray}
	\label{eqn:info_prior_likelihood}
 	P( D, \MPk\vert \MCk) & = & Z \, 2^{-\beta_k \mathbf{(} \KLd [Q \| P ]
 	 + \hmu [Q]\mathbf{)}} \\
	& \times & 2^{+\Asize^{k+1} \mathbf{(} \KLd [ U \| P ] 
	+ \hmu [U]\mathbf{)}} ~, \nonumber
\end{eqnarray}
where $\beta_k = \sum_{\hk,s} \left[ \nsks + \asks \right]$ and the
distribution of true parameters is
$P = \{ \psk, \psks \}$. The distributions $Q$ and $U$ are given by
\begin{eqnarray}
	\label{eqn:pme_distribution}
	Q & = & \left\{ \qsk = \frac{\nsk+\ask}{\beta_k} , \right. \\ 
	  & &	\left. \qsks = \frac{\nsks + \asks}{\nsk + \ask} \right\}
	  \nonumber \\
	\label{eqn:uniform_distribution}
	U & = & \left\{ \usk = \frac{1}{\Asize^k}, \usks = \frac{1}{\Asize} \right\}
	~,
\end{eqnarray}
where $Q$ is the distribution defined by the posterior mean and $U$ is a uniform
distribution. The information-theoretic quantities used above are given by
\begin{eqnarray}
	\KLd [ Q \| P ] & = & \sum_{s, \hk} \qsk \qsks \log_2 \frac{\qsks}{\psks}
	\label{eqn:conditional_KL_div} \\
	\hmu [ Q ] 	& = & - \sum_{s, \hk} \qsk \qsks \log_2 \qsks ~. 
	\label{eqn:entropy_rate_estimate}
\end{eqnarray}
The form of~\eqnref{eqn:info_prior_likelihood} and its relation to the evidence 
suggests a connection to statistical mechanics: The evidence 
$P(D \vert \MCk) = \int d\MPk P( D, \MPk\vert \MCk)$ is a partition function 
$\mathcal{Z} = P( D \vert \MCk)$.  Using conventional techniques, the
expectation and variance of the ``energy''
\begin{equation}
\label{eqn:info_energy}
\EQP = \KLd [Q \| P ] + \hmu [Q]
\end{equation}
are obtained by taking derivatives of the logarithm of the partition function
with respect to $\beta_k$:
\begin{eqnarray}
	\avg{\, \EQP \, }{\rm{post}} 
	& = &
	- \frac{1}{\log 2}
	\frac{\partial}{\partial \beta_k} \, \log \mathcal{Z}
	\label{eqn:info_mean_energy}\\
 	\var{\, \EQP \, }{\rm{post}}
	& = &
	\frac{1}{\log 2}
 	\frac{\partial^2}{\partial \beta_k^2} \, \log \mathcal{Z}
	~. 
	\label{eqn:info_variance_energy}
\end{eqnarray}
The factors of $\log 2$ in the above expressions come from the decision to use
base 2 logarithms in the definition of our information-theoretic quantities. 
This results in values in \emph{bits} rather than \emph{nats}~\cite{Cover1991}. 

To evaluate the above expression, we take advantage of the known form for the 
evidence provided in~\eqnref{eqn:evidence}.  With the definitions $\alpha_k =
\sum_{\hk} \ask$ and
\begin{equation}
	\label{eqn:prior_distribution}
	R = \left\{ \rsk = \frac{\ask}{\alpha_k} , 
	\rsks = \frac{\asks}{\ask} \right\}
\end{equation}
the negative logarithm of the partition function can be written
\begin{eqnarray}
	- \log \mathcal{Z} & = & \sum_{\hk,s} \log \Gamma 
	\left[ \alpha_k \rsk \rsks \right]
	\\ & - & \sum_{\hk} \log \Gamma \left[ \alpha_k \rsk \right] 
	+  \sum_{\hk} \log \Gamma \left[ \beta_k \qsk \right] \nonumber \\
	& - & \sum_{\hk,s} \log \Gamma 
	\left[ \beta_k \qsk \qsks \right]. \nonumber
\end{eqnarray}

From this expression, the desired expectation is found by taking derivatives
with respect to $\beta_k$; we find that
\begin{gather}
	\avg{\, \EQP \, }{\rm{post}} 
		= \frac{1}{\log 2}
		\sum_{\hk} \qsk \psi^{(0)} \left[ \beta_k \qsk \right] 
		\nonumber \\
	-  \frac{1}{\log 2} \sum_{\hk,s} \qsk \qsks \psi^{(0)} 
	\left[ \beta_k \qsk \qsks \right]~. \nonumber \\
	\label{eqn:average_info}
\end{gather} 
The variance is obtained by taking a second derivative with respect to
$\beta_k$, producing

\begin{gather}
	\var{\, \EQP \, }{\rm{post}}  = 
	- \frac{1}{\log 2} \sum_{\hk} \qsk^2 \psi^{(1)} \left[ \beta_k \qsk \right]
	\nonumber \\
	+  \frac{1}{\log 2} \sum_{\hk,s} \qsk^2 \qsks^2 \psi^{(1)} 
	\left[ \beta_k \qsk \qsks \right]. \nonumber \\
	\label{eqn:variance_info}
\end{gather} 
In both of the above the polygamma function is defined $\psi^{(n)}(x) = 
d^{n+1}/dx^{n+1} \log \Gamma(x)$. (For further details, consult a reference 
such as~\cite{Abramowitz1965}.)

From the form of~\eqnref{eqn:average_info} 
and~\eqnref{eqn:variance_info}, the meaning is not immediately clear. We can 
use an expansion of the $n=0$ polygamma function
\begin{equation}
\psi^{(0)}(x) = \log x - 1/2x + \mathcal{O}(x^{-2}) ~,
\end{equation}
valid for $x \gg 1$, however, to obtain an asymptotic form
for~\eqnref{eqn:average_info}; we find
\begin{gather}
 	\avg{\, \EQP \, }{\rm{post}} = 
 	H[ \qsk \qsks ] - H[\qsk] \nonumber \\
	+ \frac{1}{2\beta_k} \Asize^k(\Asize -1 )
	+ \mathcal{O}(1/ \beta_k^2)
	\label{eqn:average_info_asymptotic}.
\end{gather}
From this we see that the first two terms make up the entropy 
rate $\hmu [ Q ] = H[ \qsk \qsks ] - H[\qsk]$ and the last 
term is associated with the conditional relative entropy between the posterior
mean distribution $Q$ and true distribution $P$.

In summary, we have found the average of conditional relative entropy and
entropy rate with respect to the posterior density.  This was accomplished by
making connections to statistical mechanics through type theory.  Unlike
sampling from the posterior to estimate the entropy rate, this method results
in an analytic form which approaches $\hmu [ P ]$ as the inverse of the data
size. This method for approximating $\hmu$ also provides a computational
benefit. No eigenstates have to be found from the Markov transition matrix,
allowing for the storage of values in sparse data structures. This provides
a distinct computational advantage when large orders or alphabets are
considered.  

Finally, it might seem awkward to use the expectation 
of~\eqnref{eqn:info_energy} for estimation of the entropy rate.  This method
was chosen because it is the form that naturally appears in writing down the
likelihood-prior combination in~\eqnref{eqn:info_prior_likelihood}.  As a result
of using this method, most of the results obtained above are without
approximation.  We were also able to show this expectation converges to the
desired value in a well behaved manor.

\vspace{-0.125in}
\section{Examples}
\vspace{-0.125in}

To explore how the above produces a robust inference procedure, let's now
consider the statistical inference of a series of increasingly complex data
sources. The first, called the \emph{golden mean} process, is a first-order
Markov chain. The second data source is called the \emph{even process} and
cannot be represented by a Markov chain with finite order. However, this source
is a deterministic HMM, meaning that the current state and next output symbol
uniquely determine the next state.  Finally, we consider the \emph{simple
nondeterministic source}, so named since its smallest representation is as
a nondeterministic HMM. (Nondeterminism here refers to the HMM structure: the
current state and next output symbol do not uniquely determine the next state. 
This source is represented by an infinite-state deterministic HMM
\cite{Crutchfield1994,Upper1997}.)

The golden mean, even, and simple nondeterministic processes can all be written
down as models with two internal states---call them $A$ and $B$.  However, the
complexity of the data generated from each source is of markedly different
character. Our goal in this section is to consider the three main steps in
inference to analyze them. First, we consider inference of a first-order Markov
chain to demonstrate the
estimation of model parameters with uncertainty.  Second, we consider model
comparison for a range of orders $k$.  This allows us to discover structure in
the data source even though the true model class cannot be captured in all 
cases. Finally, we consider estimation of entropy rates from these data sources,
investigating how randomness is expressed in them.

While investigating these processes we consider average data counts,
rather than sample counts from specific realizations, as we want
to focus specifically on the average performance of Bayesian inference.  To
do this we take advantage of the known form of the sources. Each is described
by a transition matrix $T$, which gives transitions between states
$A$ and $B$:
\begin{equation}
	\label{eqn:transition_matrix_definition}
	T = \left[ \begin{array}{cc}
	p(A\vert A) & p(B\vert A) \\
	p(A\vert B) & p(B\vert B) 
	\end{array}
	\right] \;.
\end{equation}
Although two of our data sources are not finite Markov chains, the transition
matrix between internal states is Markov.  This means the matrix
is \emph{stochastic} (all rows sum to one) and we are guaranteed an eigenstate
$\vec{\pi}$ with eigenvalue one: $\vec{\pi} \, T = \vec{\pi}$.  This eigenstate
describes the asymptotic distribution over internal states:
$\vec{\pi} = \left[ p(A), p(B) \right]$.

The transition matrix can be divided into labeled matrices $T^{(s)}$ which
contain those elements of $T$ that output symbol $s$. For our binary data
sources one has
\begin{equation}
	\label{eqn:transition_matrix}
	T = T^{(0)} + T^{(1)}.
\end{equation}
Using these matrices, the average probability of words can be estimated for
each process of interest. For example, the probability of word $01$ can be
found using
\begin{equation}
p(01) = \vec{\pi} \, T^{(0)}T^{(1)} \vec{\eta} ~,
\end{equation}
where $\vec{\eta}$ is a column vector with all $1$'s. In this way, for any
data size $N$, we estimate the average count for a word as
\begin{equation}
\nsks = (N-k)~p(\hk s) ~.
\end{equation}
Average counts, obtained this way, will be the basis for all of
the examples to follow.

In the estimation of the true entropy rate for the examples we use the formula
\begin{equation}
	h_{\mu} = - \sum_{v \in \{A,B\}} p(v)
	\sum_{s \in \mathcal{A}} ~p(s\vert v) \log_2 p(s\vert v)
	\label{eqn:entropy_rate}
\end{equation}
for the the golden mean and even processes, where
$p(s\vert v) = T^{(s)}_{v \cdot}$ is the probability of a letter $s$ given the
state $v$ and $p(v)$ is the asymptotic probability of the state $v$ which can be
found as noted above. For the simple nondeterministic source this closed-form
expression cannot be applied and the entropy rate must be found using more
involved methods; see~\cite{Crutchfield1994} for further details.

\subsection{Golden mean process: In-class modeling}

\begin{figure}[htb]
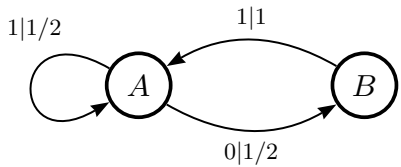

\begin{center}
			\SetStateLabelScale{1.6}
			\SetStateLineWidth{1.4pt}
			\SetEdgeLabelScale{1.4}
			\SetEdgeLineWidth{0.75pt}
		
		\begin{VCPicture}{(0,0)(5,2)}
			\ChgStateLabelScale{0.8}
				\State[A]{(1,0)}{A}
				\State[B]{(4,0)}{B}
			\ChgEdgeLabelScale{0.7}
				\LoopW{A}{ 1 | 1/2 } 
				\LArcR[0.5]{B}{A}{ 1 | 1 } 
				\LArcR[0.5]{A}{B}{ 0 | 1/2 }
		\end{VCPicture}
\end{center}
\vspace{0.5in}
\caption{A deterministic hidden Markov chain for the golden mean process.
  Edges are labeled with the output symbol and the transition probability:
  \emph{symbol} $\vert$ \emph{probability}.
  }
\label{fig:golden_mean}
\end{figure}

The \emph{golden mean process} can be represented by a simple $1$st-order
Markov chain over a binary alphabet characterized by a single (shortest)
forbidden word $s^2 = 00$. The defining labeled transition matrices for this data
source are given by
\begin{equation}
	\label{eqn:label_transition_matrix_golden_mean}
	T^{(0)} = \left[ \begin{array}{cc}
	0 & 1/2 \\
	0 & 0 
	\end{array}
	\right] \; , \;
	T^{(1)} = \left[ \begin{array}{cc}
	1/2 & 0 \\
	1	& 0 
	\end{array}
	\right]	 ~.
\end{equation}
\begfigref{fig:golden_mean} provides a graphical representation of the
corresponding hidden Markov chain. Inspection reveals a simple relation
between the \text{internal states} $A$ and $B$ and the output symbols
$0$ and $1$. An observation of $0$ indicates a transition to internal
state $B$ and a $1$ corresponds to state $A$, making this process a Markov
chain over $0$s and $1$s. 

For the golden mean the eigenstate is $\vec{\pi} = \left[ p(A), p(B)
\right] = \left( 2/3 , 1/3 \right)$.  With this vector and the labeled 
transition matrices any desired word count can be found as discussed above.

\vspace{-0.125in}
\subsubsection{Estimation of $M_1$ Parameters}
\vspace{-0.125in}

\begin{figure}[htbp]
	\centering
	\includegraphics[width=0.98\columnwidth]{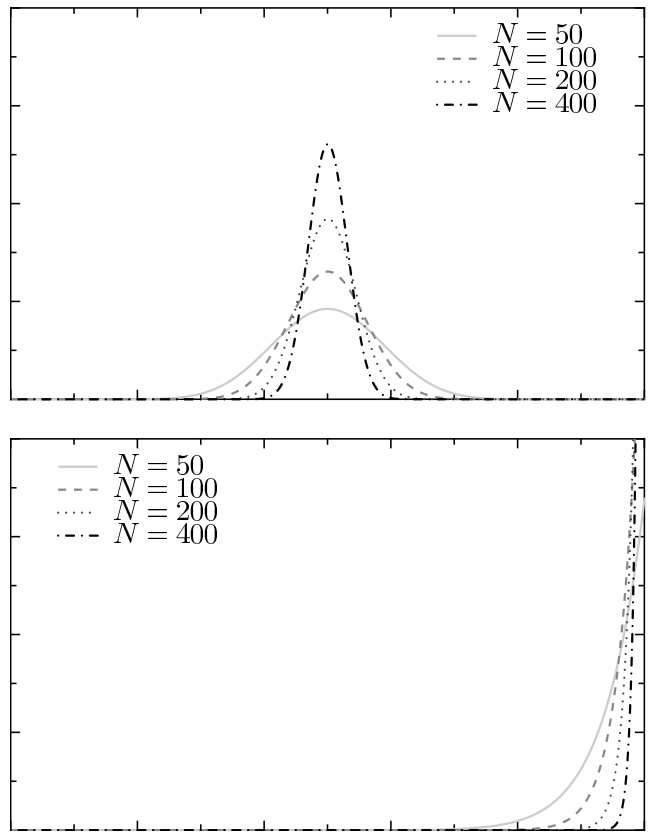}
	\caption{A plot of the inference of $M_1$ model parameters for the 
	golden mean process.  For each data sample size $N$, the marginal posterior is
	plotted for the parameters of interest: $p(0\vert 1)$ in the top panel and
	$p(1\vert 0)$ in the lower panel.  The \textit{true} values of the parameters
	are $p(0\vert 1)=1/2$ and $p(1\vert 0) = 1$.
	\label{fig:GoldenMean_ParameterEstimates}}
\end{figure}

To demonstrate the effective inference of the Markov chain parameters for the 
golden mean process we consider average counts for a variety of data sizes 
$N$.  For each size, the marginal posterior for the parameters $p(0\vert 1)$ and
$p(1\vert 0)$ is plotted in~\figref{fig:GoldenMean_ParameterEstimates}.  The
results demonstrate that the shape of the posterior effectively
describes the distribution of possible model parameters at each $N$ and converges
to the correct values of $p(0\vert 1)=1/2$ and $p(1\vert 0)=1$ with increasing
data.

Point estimates with a variance can be provided for each of the parameters, but
these numbers by themselves can be misleading.  However, the estimate obtained
by using the mean and variance of the posterior are a more effective description
of the inference process than a maximum likelihood estimate with estimated
error given by a Gaussian approximation of the likelihood alone.
As~\figref{fig:GoldenMean_ParameterEstimates} demonstrates, in
fact, a Gaussian
approximation of uncertainty is an ineffective description of our knowledge
when the Markov chain parameters are near their upper or lower limits at $0$
and $1$. Probably the most effective set of numbers to provide consists of the
mean of the posterior and a region of confidence. These would most accurately
describe asymmetries in the uncertainty of model parameters. Although we will
not do that here, a brief description of finding regions of confidence is
provided in~\appref{app:dirichlet}.

\vspace{-0.125in}
\subsubsection{Selecting the Model Order $k$}
\vspace{-0.125in}

Now consider the selection of the appropriate order $k$ from golden mean
realizations.  As discussed above, the golden mean process is a first order
Markov chain with $k=1$.  As a result, we would expect model comparison to
select this order from the available possibilities. To demonstrate this,
we consider orders $k=1-4$ and perform model comparison with a uniform prior
over orders (\eqnref{eqn:best_model_uniform_prior}) and with a penalty for the 
number of model parameters (\eqnref{eqn:best_model_df_penalty_prior}).

\begin{figure}[htbp]
	\centering
	\includegraphics[width=0.98\columnwidth]{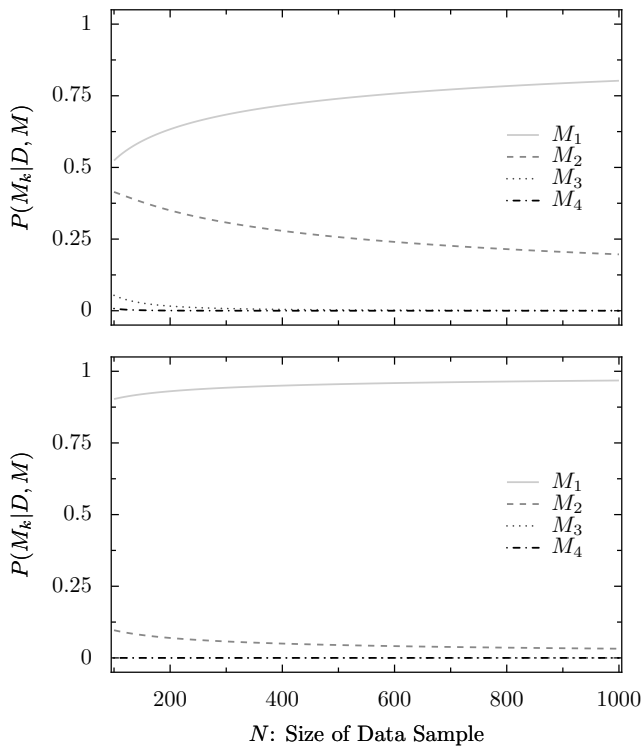}
	\caption{Model comparison for Markov chains of order $k=1-4$ using 
	average counts from the golden mean process.  Sample sizes from $N=100$ to 
	$N=1,000$ in steps of $\Delta N=5$ are used to generate these plots.  The top panel
	displays the model probabilities using a uniform prior over orders $k$.  The
	bottom panel displays the effect of a penalty for model size.
	\label{fig:GoldenMean_ModelComparison}}
\end{figure}
 
The results of the model comparisons are given 
in~\figref{fig:GoldenMean_ModelComparison}.  The top panel shows the probability
for each order $k$ as a function of the sample size, using a uniform prior.  For
this prior over orders, $M_1$ is selected with any reasonable amount of
data.   However, there does seem to be a possibility to over-fit for small data
size $N \leq 100$.  The bottom panel shows the model probability with a penalty
prior over model order $k$.  This removes the over-fitting at small data sizes
and produces an offset which must be overcome by the data before higher $k$ is
selected.  This example is not meant to argue for the penalty prior over model
orders.  In fact, Bayesian model comparison with a uniform prior does an
effective job using a relatively small sample size.

\vspace{-0.125in}
\subsubsection{Estimation of Entropy Rate}
\vspace{-0.125in}

We can also demonstrate the convergence of  the average for
$E(Q,P)=D[ Q \| P ] + \hmu [Q]$ given in~\eqnref{eqn:average_info} to the
correct entropy rate for the golden mean process.  We choose to show this
convergence for all orders $k=1-4$ discussed in the previous section.  This
exercise demonstrates that all orders greater than or equal to $k = 1$
effectively capture the entropy rate. However, the convergence to the correct
values for higher-order $k$ takes more data because of a larger initial value of
$D[ Q \| P ]$.  This larger value is simply due to the larger number of
parameters for higher-order Markov chains.

\begin{figure}[htbp]
	\centering
	\includegraphics[width=0.98\columnwidth]{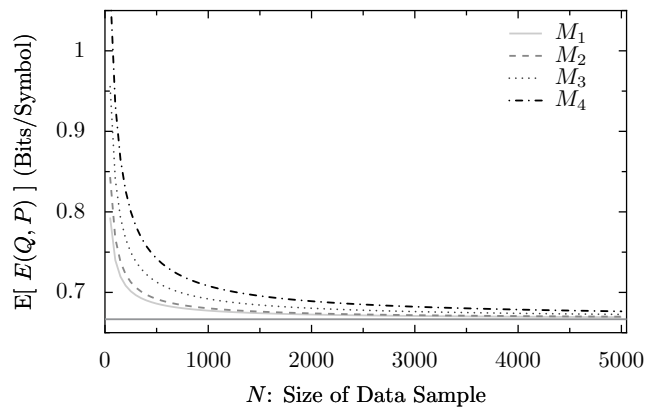}
	\caption{The convergence of $\avg{\, E(Q,P) \, }{\rm{post}}$ to the true 
	entropy rate $\hmu = 2/3$ bits per symbol (indicated by the gray horizontal 
	line) for the the golden mean process.  As demonstrated
	in~\eqnref{eqn:average_info_asymptotic}, the conditional relative 
	entropy $D[Q \| P ] \rightarrow 0$ as $1/N$.  This results in 
	the convergence of $\hmu [Q]$ to the true entropy rate.
	\label{fig:GoldenMean_InfoTheory}}
\end{figure}

In evaluating the value of $D[Q \| P ] + \hmu [Q]$ for different sample lengths,
we expect that the \pme \, estimated $Q$ will converge to the true distribution
$P$.  As a result, the conditional relative entropy should go to zero with
increasing $N$.  For the golden mean process, the known value of the entropy 
rate is $\hmu =2/3$ bits per symbol.  Inspection 
of~\figref{fig:GoldenMean_InfoTheory} demonstrates the expected convergence of the 
average from~\eqnref{eqn:average_info} to the true entropy rate.

The result of our model comparison from the previous section could also be used
in the estimation of the entropy rate.  As we saw 
in~\figref{fig:GoldenMean_ModelComparison}, there are ranges of sample length $N$
where the probability of orders $k=1,2$ are both nonzero.  In principle, an
estimate of $\hmu$ should be made by weighting the values obtained for each
$k$ by the corresponding order probability $P(\MCk \vert D, \mathcal{M})$.  As
we can see from~\figref{fig:GoldenMean_InfoTheory}, the estimates of the entropy
rate for $k=1,2$ are also very similar in this range of $N$.  As a result, this
additional step would not have a large effect for entropy rate estimation.

\subsection{Even process: Out-of-class modeling}

We now consider a more difficult data source called the \emph{even process}.
The defining labeled transition matrices are given by
\begin{equation}
	\label{eqn:label_transition_matrix_even}
	T^{(0)} = \left[ \begin{array}{cc}
	1/2 & 0 \\
	0	& 0 
	\end{array}
	\right] \; , \;
	T^{(1)} = \left[ \begin{array}{cc}
	0 & 1/2 \\
	1	& 0 
	\end{array}
	\right]~.
\end{equation}

\begin{figure}[htb]
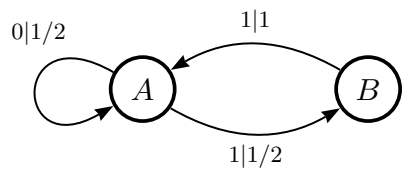

\begin{center}
			\SetStateLabelScale{1.6}
			\SetStateLineWidth{1.4pt}
			\SetEdgeLabelScale{1.4}
			\SetEdgeLineWidth{0.75pt}
		
		\begin{VCPicture}{(0,0)(5,2)}
			\ChgStateLabelScale{0.8}
				\State[A]{(1,0)}{A}
				\State[B]{(4,0)}{B}
			\ChgEdgeLabelScale{0.7}
				\LoopW{A}{ 0 | 1/2 } 
				\LArcR[0.5]{B}{A}{ 1 | 1 } 
				\LArcR[0.5]{A}{B}{ 1 | 1/2 }
		\end{VCPicture}
\end{center}
\vspace{0.5in}
\caption{Deterministic hidden Markov chain representation of the even process.
  This process cannot be represented as a finite-order (nonhidden) Markov chain
  over the output symbols $0$s and $1$s. The set of irreducible forbidden words
  $\mathcal{F} = \{01^{2n+1}0: n=0,1,2,\ldots \}$ reflects the fact that the
  process generates blocks of $1$'s, bounded by $0$s, that are \emph{even} in
  length, at any length.  
\label{fig:even}}
\end{figure}

As can be seen in~\figref{fig:even}, the node-edge structure is identical to
the golden mean process but the output symbols on the edges have been changed
slightly.  As a result of this shuffle, the states $A$ and $B$ can no longer be
associated with a simple sequence of $0$'s and $1$'s.  Whereas the golden mean
has the irreducible set of forbidden words $\mathcal{F} = \{00\}$, the even
process has a countably infinite set
$\mathcal{F} = \{01^{2n+1}0: n=0,1,2,\ldots \}$
\cite{Crutchfield1994}.

In simple terms, the even process produces blocks of $1$'s which are even in
length. This is a much more complicated type of memory than we saw in
the golden mean process.  For the Markov chain model class, where a word of
length $k$ is used to predict the next letter, this would require an
infinite-order $k$. It would be necessary to keep track of all even and odd
strings of $1$'s, irrespective of the length. As a result, the properties of
the even process mean that a finite Markov chain \textit{cannot} represent
this data source.  

This example is then a demonstration of what can be learned in a case of
out-of-class modeling. We are interested, therefore, in how well Markov
chains approximate the even process. We
expect that model comparison will select larger $k$ as the size of the data 
sample increases.  Does the model selection tells us anything about the 
underlying data source despite the inability to exactly capture its properties?
As we will see, we do obtain intriguing hints of the true nature of the even
process from model comparison.  Finally, can we estimate the entropy rate of
the process with a Markov chain?  As we will see, a high $k$ is needed to do
this effectively.

\subsubsection{Estimation of $M_1$ Parameters}

In this section we consider an $M_1$ approximation of the even process.
We expect the resulting model to accurately capture length-$2$ word 
probabilities as $N$ increases.  In this example, we consider the \emph{true}
model to be the best approximation possible by a $k=1$ Markov chain.  From the
labeled transition matrices given above we can calculate the appropriate
values for $p(0\vert 1)$ and $p(1\vert 0)$ using the methods described above.  
Starting from the asymptotic distribution $\vec{\pi} = \left[ p(A), p(B)\right]
= \left[ 2/3, 1/3 \right]$ we obtain $p(0\vert 1)=p(10)/p(1)=1/4$ and $p(1\vert
0)=p(01)/p(0)=1/2$. 

\begin{figure}[htbp]
	\centering
	\includegraphics[width=0.98\columnwidth]{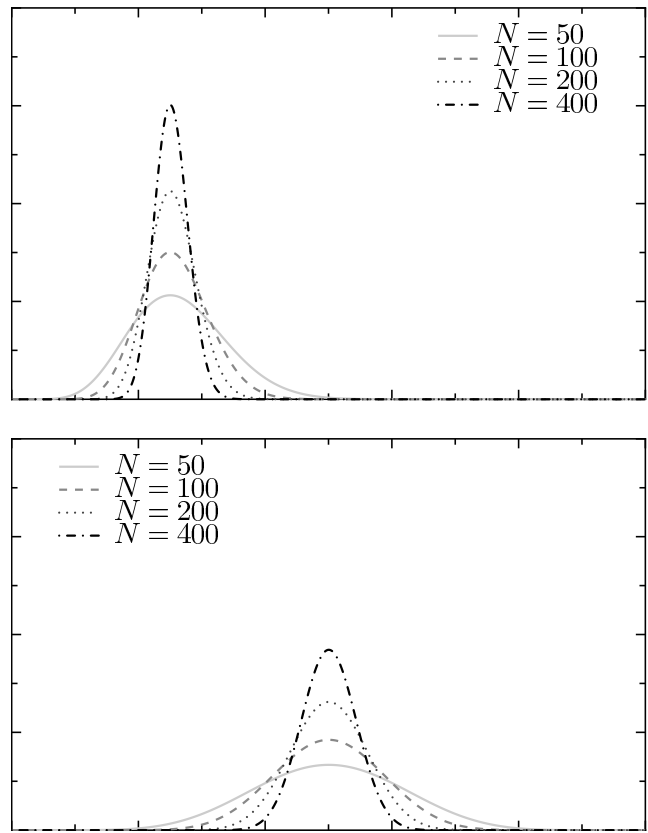}
	\caption{A plot of the inference of $M_1$ model parameters for the even 
	process.  For a variety of sample sizes $N$, the marginal posterior for
	$p(0\vert 1)$ (top panel) and $p(1\vert 0)$ (bottom panel) are shown.  The
	\textit{true} values of the parameters are $p(0\vert 1)=1/4$ and 
	$p(1\vert 0) = 1/2$.
	\label{fig:Even_ParameterEstimates}}
\end{figure}

As we can see from~\figref{fig:Even_ParameterEstimates}, a first-order Markov
chain can be inferred without difficulty.  The values obtained are exactly as
expected.  However, these values do not tell us much about the nature
of the data source by themselves. This points to the important role of model
comparison and entropy rate estimation in understanding this data.

\subsubsection{Selecting the Model Order $k$}

Now consider the selection of Markov chain order $k=1-4$ for a range of data
sizes $N$. Recall that the even process cannot be represented by a finite-order
Markov chain over the output symbols $0$ and $1$. As a consequence, we expect
higher $k$ to be selected with increasing data $N$, as more data statistically
justifies more complex models. This is what happens, in fact, but the way in
which orders are selected as we increase $N$ provides structural information
we could not obtain from the inference of a Markov chain of fixed order.

\begin{figure}[htbp]
	\centering
	\includegraphics[width=0.98\columnwidth]{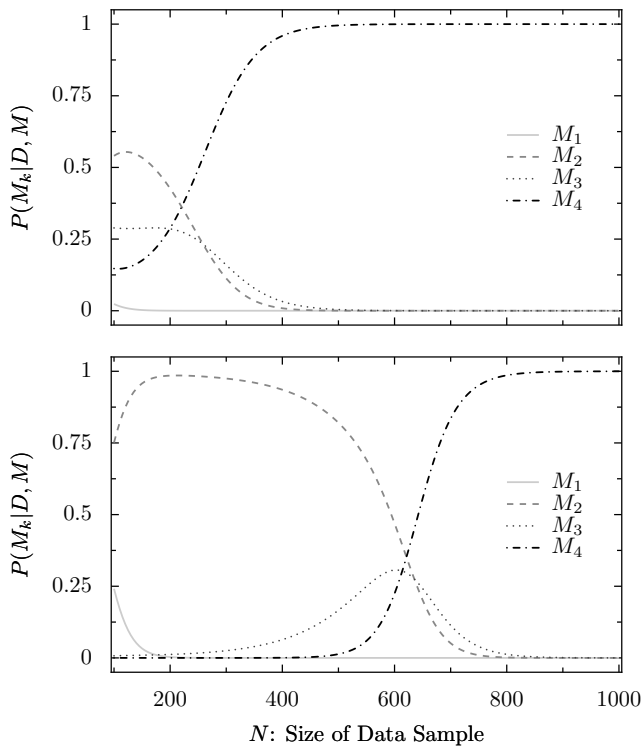}
	\caption{Model comparison for Markov chains of order $k=1-4$ for 
	average data from the even process.  The top panel shows the model 
	comparison with a uniform prior over the possible orders $k$.  The bottom 
	panel demonstrates model comparison with a penalty for the number of model
	parameters.  In both cases the $k=4$ model is chosen over lower orders as the
	amount of data available increases.
	\label{fig:Even_ModelComparison}}
\end{figure}

If we consider~\figref{fig:Even_ModelComparison}, an interesting pattern becomes
apparent.  Orders with even $k$ are preferred over odd. In this way model
selection is hinting at the underlying structure of the source. The Markov
chain model class cannot represent the even process in a compact way, but
inference and model comparison combined provide useful information about
the hidden structure of the source.

In this example we also have regions where the probability of multiple orders $k$
are equally probable.  The sample size at which this occurs depends on the prior
over orders which is employed.  When this happens, properties estimated from the
Markov chain model class should use a weighted sum of the various orders. As we
will see in the estimation of entropy rates, this is not as critical. At
sample sizes where the order probabilities are similar, the estimated entropy
rates are also similar.

\subsubsection{Estimation of Entropy Rate}

\begin{figure}[htbp]
	\centering
	\includegraphics[width=0.98\columnwidth]{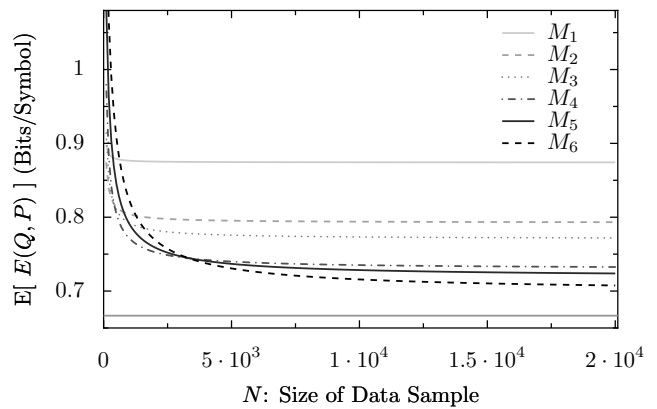}
	\caption{The convergence of $\avg{\, D[Q \| P ] + \hmu [Q] \, }{\rm{post}}$ 
	to the true entropy rate $\hmu = 2/3$ bits per symbol for the the even
	process.  The true value is indicated by the horizontal gray line.  Experience
	shows that a $k=10$ Markov chain is needed to effectively approximate the true
	value of $\hmu$.
	\label{fig:Even_InfoTheory}}
\end{figure}

Entropy rate estimation for the even process turns out to be a more
difficult task than one might expect.  In~\figref{fig:Even_InfoTheory} we see
that Markov chains of orders $1-6$ are unable to effectively capture the true
entropy rate.  In fact, experience shows that an order $k=10$ Markov chain or
higher is needed to get close to the true value of $\hmu = 2/3$ bits per symbol.
Note also the factor of $20$ longer realizations that are required compared,
say, to the golden mean example.

As discussed above, a weighted sum of 
$\avg{\, D[Q \| P ] + \hmu [Q] \, }{\rm{post}}$ could be employed in this
example.  For the estimate this is not critical because the different orders
provide roughly the same value at these points.  In fact, these points
correspond to where the estimates of $E(Q,P)$ cross 
in~\figref{fig:Even_InfoTheory}. They are samples sizes where apparent
randomness can be explained by structure and increased order $k$.

\subsection{Simple Nondeterministic Source: Out-of-class modeling}

\begin{figure}[htb]
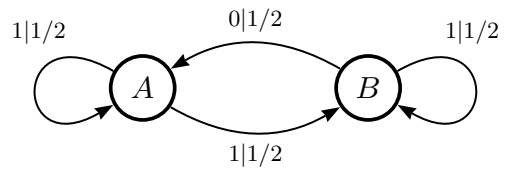

\begin{center}
			\SetStateLabelScale{1.6}
			\SetStateLineWidth{1.4pt}
			\SetEdgeLabelScale{1.4}
			\SetEdgeLineWidth{0.75pt}
		
		\begin{VCPicture}{(0,0)(5,2)}
			\ChgStateLabelScale{0.8}
				\State[A]{(1,0)}{A}
				\State[B]{(4,0)}{B}
			\ChgEdgeLabelScale{0.7}
				\LoopW{A}{ 1 | 1/2 }
				\LoopE{B}{ 1 | 1/2 } 
				\LArcR[0.5]{B}{A}{ 0 | 1/2 } 
				\LArcR[0.5]{A}{B}{ 1 | 1/2 }
		\end{VCPicture}
\end{center}
\vspace{0.5in}
\caption{A hidden Markov chain representation of the simple nondeterministic
  process. This example also cannot be represented as a finite-order Markov
  chain over outputs $0$ and $1$. It, however, is more complicated than the
  two previous examples: Only the observation of a $0$ provides the observer
  with information regarding the internal state of the underlying process;
  observing a $1$ leaves the internal state ambiguous.
\label{fig:sns}}
\end{figure}

The simple nondeterministic source adds another level of challenge to inference.
As its name suggests, it is described by a nondeterministic HMM.  
Considering~\figref{fig:sns} we can see that a $1$ is produced on every
transition except for the $B \rightarrow A$ edge.  This means there are many
paths through the internal states that produce the same observable sequence of
$0$s and $1$s. The defining labeled transition matrices for this process are
given by
\begin{equation}
	\label{eqn:label_transition_matrix_sns}
	T^{(0)} = \left[ \begin{array}{cc}
	0 & 0 \\
	1/2	& 0 
	\end{array}
	\right] \; , \;
	T^{(1)} = \left[ \begin{array}{cc}
	1/2 & 1/2 \\
	0	& 1/2 
	\end{array}
	\right]~.
\end{equation}

Using the state-to-state transition matrix $T=T^{(0)}+T^{(1)}$, we find the
asymptotic distribution for the hidden states to be
$\vec{\pi} = \left[ p(A), p(B) \right] = \left[1/2, 1/2 \right]$. Each of
the hidden states is equally likely; however, a $1$ is always produced from
state $A$, while there is an equal chance of obtaining a $0$
or $1$ from state $B$.

\subsubsection{Estimation of $M_1$ Parameters}

\begin{figure}[htbp]
	\centering
	\includegraphics[width=0.98\columnwidth]{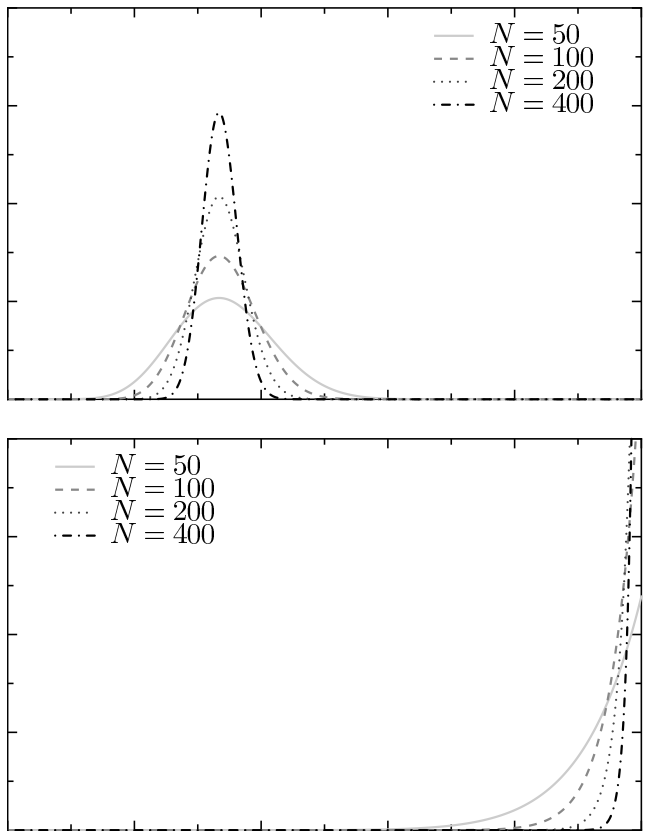}
	\caption{Marginal density for $M_1$ model parameters for the
	simple nondeterministic process:  The curves for each data size $N$
	demonstrate a well behaved convergence to the correct values:
	$p(0\vert 1)=1/3$ and $p(1\vert 0) = 1$.
	\label{fig:SNS_ParameterEstimates}}
\end{figure}

Using the asymptotic distribution derived above, the parameters of an inferred
first-order Markov chain should approach $p(0\vert 1)=p(10)/p(1)=1/3$ and
$p(1\vert 0)=p(01)/p(0)=1$.  As we can see
from~\figref{fig:SNS_ParameterEstimates}, the inference
process captures these values very effectively despite the out-of-class data
source.  

\subsubsection{Selecting the Model Order $k$}

\begin{figure}[htbp]
	\centering
	\includegraphics[width=0.98\columnwidth]{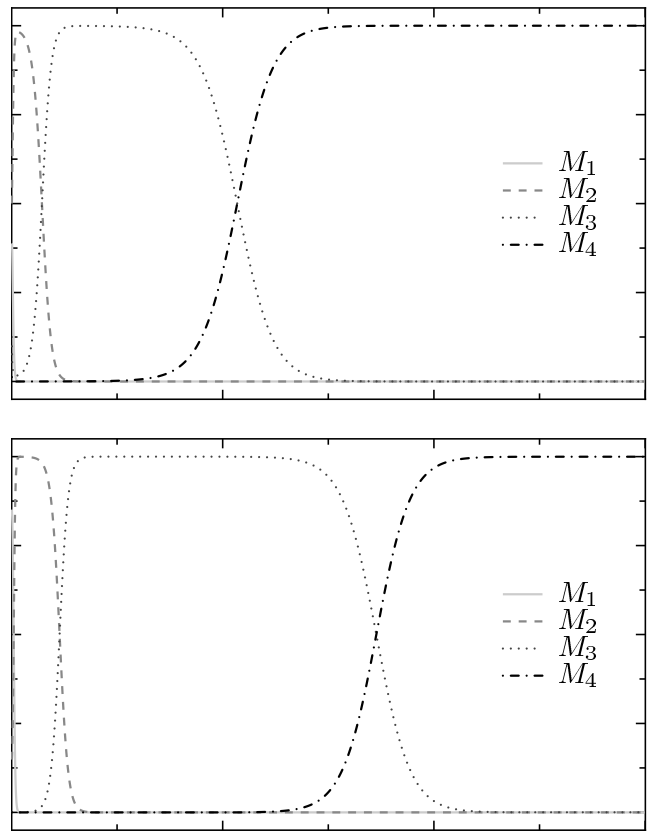}
	\caption{Model comparison for Markov chains of order $k=1-4$ for 
	data from the simple nondeterministic process.  The top panel 
	shows the model comparison with a uniform prior over the possible orders 
	$k$.  The bottom panel demonstrates model comparison with a penalty for the 
	number of model parameters.  Note the scale on the horizontal axis---it
	takes much more data for the model comparison to pick out higher orders
	for this process compared to the previous examples.
	\label{fig:SNS_ModelComparison}}
\end{figure}

Here we consider the comparison of Markov chain models of orders $k=1-4$ when
applied to data from the simple nondeterministic source.  As with the even
process, we expect increasing order to be selected as the amount of available
data increases.  In~\figref{fig:SNS_ModelComparison} we see that this is
exactly what happens.

Unlike the even process, there is no preference for even orders.  Instead, we
observe a systematic increase in order with larger data sets.  We do note that
the amount of data need to select a higher order does seem to be larger than for
the even process.  Here the distribution over words is more important and more
subtle than the support of the distribution (those words with positive
probability).

\subsubsection{Estimation of Entropy Rate}

Estimation of the entropy rate for the simple nondeterministic source provides
an interesting contrast to the previous examples. As discussed when introducing
the examples, this data source is a nondeterministic HMM and the entropy rate
cannot be directly calculated using~\eqnref{eqn:entropy_rate}
\cite{Blackwell1957}. However, a
value of $\hmu \approx 0.677867$ bits per symbol has been obtained
in~\cite{Crutchfield1994}. 

\begin{figure}[htbp]
	\centering
	\includegraphics[width=0.98\columnwidth]{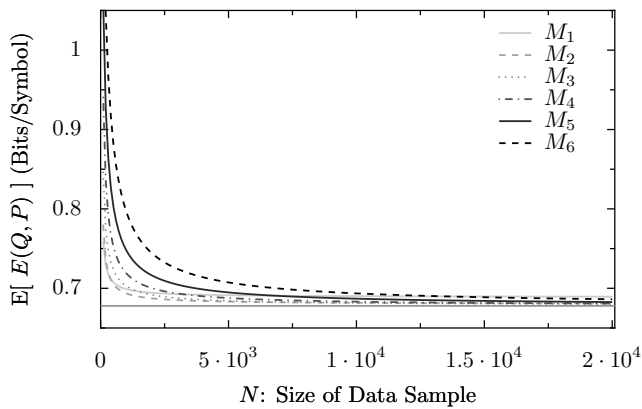}
	\caption{The convergence of $\avg{\, D[Q \| P ] + 
	\hmu [Q] \, }{\rm{post}}$ to the true entropy rate $\hmu \approx 0.677867$ 
	bits per symbol for the simple nondeterministic source.  The true value is
	indicated by the gray horizontal line.
	\label{fig:SNS_InfoTheory}}
\end{figure}

\begfigref{fig:SNS_InfoTheory} shows the results of entropy-rate estimation
using Markov chains of order $k=1-6$.  These results demonstrate that the
entropy rate can be effectively estimated with low-order $k$ and relatively
small data samples.  This is an interesting result, as we might expect
estimation of the entropy rate to be most difficult in this example.  Instead we
find that the even process was a more difficult test case.

\section{Discussion}

The examples presented above provide several interesting lessons in inference,
model comparison, and estimating randomness. The combination of these three
ideas applied to a data source provides information and intuition about the
structure of the underlying system, even when modeling out-of-class processes.

In the examples of $\MC_{1}$ estimates for each of the sources we see that
the Bayesian methods provide a powerful and consistent description of Markov
chain model parameters.  The marginal density accurately describes the
uncertainty associated with these estimates, reflecting asymmetries which point
estimation with error bars cannot capture.  In addition, methods described
in~\appref{app:dirichlet} can be used to generate regions of confidence of any
type. 

Although the estimates obtained for the Markov chain model parameters were
consistent with the data source for words up to length $k+1$, they did not capture
the true nature of the system under consideration.  This demonstrates that
estimation of model parameters without some kind of model comparison can be very
misleading.  Only with the comparison of different orders did some indication
of the true properties of the data source become clear.  Without this step, 
misguided interpretations are easily obtained.

For the golden mean process, a $k=1$ Markov chain, the results of model
comparison were predictably uninteresting.  This is a good indication that the
correct model class is being employed.  However, with the even process a much
more complicated model comparison was found.  In this case, a selection of even
$k$ over odd hinted at the distinguishing properties of the source. In a
similar way, the results of model comparison for the simple nondeterministic
source selected increasing order with larger $N$.  In both out-of-class
modeling examples, the increase in selected order without end is a good
indication that the data source is not in the Markov chain class. (A parallel
technique is found in \emph{hierarchical $\epsilon$-machine reconstruction}
\cite{Crutchfield1994}.) Alternatively, there is an indication that
very high-order dependencies are important in the description of the process. 
Either way, this information is important since it gives an indication to the
modeler that a more complicated dynamic is at work and all results must be
treated with caution.

Finally, we considered the estimation of entropy rates for the example data
sources.  In two of the cases, the golden mean process and the simple
nondeterministic source, short data streams were adequate.  This is not
unexpected for the golden mean, but for the simple nondeterministic source this
might be considered surprising.  For the even process, the estimation of the
entropy rate was markedly more difficult.  For this data source, the countably
infinite number of forbidden words makes the support of the word distribution
at a given length important.  As a result, a larger amount of data and a
higher-order Markov chain are needed to find a decent estimate of randomness
from that data source. In this way, each of the steps in Bayesian
inference allow one to separate structure from randomness. 

\section{Conclusion}

We considered Bayesian inference of $k$-th order Markov chain
models.  This included estimating model parameters for a given $k$, model
comparison between orders, and estimation of randomness in the form of entropy
rates.  In most approaches to inference, these three aspects are treated as
separate, but related endeavors.  However, we find them to be intimately
related.  An estimate of model parameters without a sense of whether the
correct model is being used is misguided at best.  Model comparison
provides a window into this problem by comparing various orders $k$ within the
model class.  Finally, estimating randomness in the form of an entropy rate
provides more information about the trade-off between structure and randomness. 
To do this we developed a connection to the statistical mechanical partition
function, from which averages and variances were directly calculable. For the
even process, structure was perceived as randomness and for the simple
nondeterministic source
randomness was easily estimated and structure was more difficult to find.
These insights, despite the out-of-class data, demonstrate the power of
combining these three methods into one effective tool for investigating
structure and randomness in finite strings of discrete data.

%
%
\section*{Acknowledgments}
This work was partially supported at the Center for Computational Science
and Engineering at the University of California at Davis by Intel
Corporation. Work at the Santa Fe Institute was supported under its
Computation, Dynamics, and Inference Program core grants from the
National Science and MacArthur Foundations. C.S. and A.H. acknowledge
support by the National Science Foundation Grant DMS 03-25939 ITR.

%
%
\appendix

%
%
\section{}
\label{app:Dirichlet}

\subsection{Dirichlet Distribution\label{app:dirichlet}}

We supply a brief overview of the Dirichlet distribution for completeness.  For 
more information, a reference such as~\cite{Wilks1962} should be consulted.  In 
simple terms, the Dirichlet distribution is the multinomial generalization of 
the Beta distribution.  The probability density function for $q$ elements is
given by
\begin{equation}
	\label{eqn:dirichlet_pdf}
	\text{Dir}( \{ p_{i} \} )
	=
	\frac{ \Gamma( \alpha ) }{\prod_{i=0}^{q-1} \Gamma( \alpha_{i} ) }
	\delta(1-\sum_{i=0}^{q-1} p_{i})
	\prod_{i=0}^{q-1} p_{i}^{\alpha_{i}-1}.
\end{equation}

The variates must satisfy $p_i \in [0,1]$ and $\sum_{i=0}^{q-1} p_{i} = 1$. The 
hyperparameters $\{ \alpha_{i} \}$ of the distribution, must be real and 
positive and we use the notation $\alpha = \sum_{i=0}^{q-1} \alpha_{i}$.  The
average, variance, and covariance of the parameters $p_{i}$ are
given by, respectively,
\begin{eqnarray}
	\avg{p_{j}}{} & = & \frac{ \alpha_{j} }{ \alpha }, 
	\label{eqn:dirichlet_average}\\
	\var{p_{j}}{} & = & \frac{ \alpha_{j}\left( \alpha - \alpha_{j} \right) 
	}{ \alpha^{2} \left( 1+ \alpha \right) },
	\label{eqn:dirichlet_variance}\\
	\cov{p_{j}}{p_{l}}	& = & - \frac{ \alpha_{j} \alpha_{l} 
	}{ \alpha^{2} \left( 1+ \alpha \right) } \; , \; j \neq l. 
	\label{eqn:dirichlet_covariance}
\end{eqnarray}

\subsection{Marginal distributions\label{app:dirichlet_marginal}}

An important part of understanding uncertainty in the inference process is the
ability to find regions of confidence from a marginal density.  The marginal is
obtained from the posterior by integrating out the dependence on all parameters
except for the parameter of interest.  For a Dirichlet distribution, the
marginal density is known to be a Beta distribution~\cite{Wilks1962},
\begin{equation}
	\label{eqn:beta_pdf}
	\text{Beta}( p_{i} )
	=
	\frac{ \Gamma( \alpha ) }{\Gamma( \alpha_{i} ) \Gamma( \alpha - \alpha_{i} ) }
	 p_{i}^{\alpha_{i}-1} \left( 1 - p_{i} \right)^{\alpha - \alpha_{i}-1}.
\end{equation}

\subsection{Regions of confidence from the marginal density}

From the marginal density provided in~\eqnref{eqn:beta_pdf} a cumulative
distribution function can be obtained using the incomplete Beta integral
\begin{equation}
	\Pr(p_{i} \leq x) = \int_{0}^{x} \, dp_{i} \, \text{Beta}(p_{i}) ~.
	\label{eqn:beta_cdf}
\end{equation}
Using this form, the probability that a Markov chain parameter will be between
$a$ and $b$ can be found using $\Pr( a \leq p_{i} \leq b) = \Pr( p_{i} \leq b) -
\Pr( p_{i} \leq a)$.  For a confidence level $R$, between zero and one, we then
want to find $(a,b)$ such that $R=\Pr( a \leq p_{i} \leq b)$.  The incomplete
Beta integral and its inverse can be found using computational methods,
see~\cite{Majumder1973,Majumder1973a,Cran1977,Berry1990} for details.

%
%


\begin{thebibliography}{29}
\expandafter\ifx\csname natexlab\endcsname\relax\def\natexlab#1{#1}\fi
\expandafter\ifx\csname bibnamefont\endcsname\relax
  \def\bibnamefont#1{#1}\fi
\expandafter\ifx\csname bibfnamefont\endcsname\relax
  \def\bibfnamefont#1{#1}\fi
\expandafter\ifx\csname citenamefont\endcsname\relax
  \def\citenamefont#1{#1}\fi
\expandafter\ifx\csname url\endcsname\relax
  \def\url#1{\texttt{#1}}\fi
\expandafter\ifx\csname urlprefix\endcsname\relax\def\urlprefix{URL }\fi
\providecommand{\bibinfo}[2]{#2}
\providecommand{\eprint}[2][]{\url{#2}}

\bibitem[{\citenamefont{Avery and Henderson}(1999)}]{Avery1999}
\bibinfo{author}{\bibfnamefont{P.~J.} \bibnamefont{Avery}} \bibnamefont{and}
  \bibinfo{author}{\bibfnamefont{D.~A.} \bibnamefont{Henderson}},
  \bibinfo{journal}{Appl. Stat.} \textbf{\bibinfo{volume}{48}},
  \bibinfo{pages}{53 } (\bibinfo{year}{1999}).

\bibitem[{\citenamefont{Liu and Lawrence}(1999)}]{JSLiu1999}
\bibinfo{author}{\bibfnamefont{J.~S.} \bibnamefont{Liu}} \bibnamefont{and}
  \bibinfo{author}{\bibfnamefont{C.~E.} \bibnamefont{Lawrence}},
  \bibinfo{journal}{Bioinformatics} \textbf{\bibinfo{volume}{15}},
  \bibinfo{pages}{38 } (\bibinfo{year}{1999}).

\bibitem[{\citenamefont{Crutchfield and Feldman}(1997)}]{Crutchfield1997}
\bibinfo{author}{\bibfnamefont{J.~P.} \bibnamefont{Crutchfield}}
  \bibnamefont{and} \bibinfo{author}{\bibfnamefont{D.~P.}
  \bibnamefont{Feldman}}, \bibinfo{journal}{Phys. Rev. E}
  \textbf{\bibinfo{volume}{55}}, \bibinfo{pages}{R1239 }
  (\bibinfo{year}{1997}).

\bibitem[{\citenamefont{MacKay and Peto}(1994)}]{MacKay1994}
\bibinfo{author}{\bibfnamefont{D.~J.~C.} \bibnamefont{MacKay}}
  \bibnamefont{and} \bibinfo{author}{\bibfnamefont{L.~C.~B.}
  \bibnamefont{Peto}}, \bibinfo{journal}{Nat. Lang. Eng.}
  \textbf{\bibinfo{volume}{1}} (\bibinfo{year}{1994}).

\bibitem[{\citenamefont{Crutchfield and Packard}(1983)}]{Crutchfield1983}
\bibinfo{author}{\bibfnamefont{J.~P.} \bibnamefont{Crutchfield}}
  \bibnamefont{and} \bibinfo{author}{\bibfnamefont{N.~H.}
  \bibnamefont{Packard}}, \bibinfo{journal}{Physica D}
  \textbf{\bibinfo{volume}{7D}}, \bibinfo{pages}{201 } (\bibinfo{year}{1983}).

\bibitem[{\citenamefont{Hao and Zheng}(1998)}]{BLHao1998}
\bibinfo{author}{\bibfnamefont{B.-L.} \bibnamefont{Hao}} \bibnamefont{and}
  \bibinfo{author}{\bibfnamefont{W.-M.} \bibnamefont{Zheng}},
  \emph{\bibinfo{title}{Applied Symbolic Dynamics and Chaos}}
  (\bibinfo{publisher}{World Scientific}, \bibinfo{year}{1998}).

\bibitem[{\citenamefont{Anderson and Goodman}(1957)}]{TWAnderson1957}
\bibinfo{author}{\bibfnamefont{T.~W.} \bibnamefont{Anderson}} \bibnamefont{and}
  \bibinfo{author}{\bibfnamefont{L.~A.} \bibnamefont{Goodman}},
  \bibinfo{journal}{Ann. Math. Stat.} \textbf{\bibinfo{volume}{28}},
  \bibinfo{pages}{89 } (\bibinfo{year}{1957}).

\bibitem[{\citenamefont{Billingsley}(1961)}]{Billingsley1961a}
\bibinfo{author}{\bibfnamefont{P.}~\bibnamefont{Billingsley}},
  \bibinfo{journal}{Ann. Math. Stat.} \textbf{\bibinfo{volume}{32}},
  \bibinfo{pages}{12 } (\bibinfo{year}{1961}).

\bibitem[{\citenamefont{Chatfield}(1973)}]{Chatfield1973}
\bibinfo{author}{\bibfnamefont{C.}~\bibnamefont{Chatfield}},
  \bibinfo{journal}{Appl. Stat.} \textbf{\bibinfo{volume}{22}},
  \bibinfo{pages}{7} (\bibinfo{year}{1973}).

\bibitem[{\citenamefont{Tong}(1975)}]{HTong1975}
\bibinfo{author}{\bibfnamefont{H.}~\bibnamefont{Tong}}, \bibinfo{journal}{Jour.
  Appl. Prob.} \textbf{\bibinfo{volume}{12}}, \bibinfo{pages}{488 }
  (\bibinfo{year}{1975}).

\bibitem[{\citenamefont{Katz}(1981)}]{Katz1981}
\bibinfo{author}{\bibfnamefont{R.~W.} \bibnamefont{Katz}},
  \bibinfo{journal}{Technometrics} \textbf{\bibinfo{volume}{23}},
  \bibinfo{pages}{243 } (\bibinfo{year}{1981}).

\bibitem[{\citenamefont{Rissanen}(1984)}]{JRissanen1984}
\bibinfo{author}{\bibfnamefont{J.}~\bibnamefont{Rissanen}},
  \bibinfo{journal}{IEEE Trans. Inform. Theory} \textbf{\bibinfo{volume}{30}},
  \bibinfo{pages}{629} (\bibinfo{year}{1984}).

\bibitem[{\citenamefont{Vapnik}(1999)}]{VVapnik1999}
\bibinfo{author}{\bibfnamefont{V.}~\bibnamefont{Vapnik}},
  \bibinfo{journal}{IEEE Trans. Neur. Net.} \textbf{\bibinfo{volume}{10}},
  \bibinfo{pages}{988} (\bibinfo{year}{1999}).

\bibitem[{\citenamefont{Vit{\'a}nyi and Li}(2000)}]{Vitanyi2000}
\bibinfo{author}{\bibfnamefont{P.~M.} \bibnamefont{Vit{\'a}nyi}}
  \bibnamefont{and} \bibinfo{author}{\bibfnamefont{M.}~\bibnamefont{Li}},
  \bibinfo{journal}{IEEE Trans. Inform. Theory}
  \textbf{\bibinfo{volume}{46(2)}}, \bibinfo{pages}{446}
  (\bibinfo{year}{2000}).

\bibitem[{\citenamefont{Baldi and Brunak}(2001)}]{Baldi2001}
\bibinfo{author}{\bibfnamefont{P.}~\bibnamefont{Baldi}} \bibnamefont{and}
  \bibinfo{author}{\bibfnamefont{S.}~\bibnamefont{Brunak}},
  \emph{\bibinfo{title}{Bioinformatics: The Machine Learning Approach}}
  (\bibinfo{publisher}{MIT Press}, \bibinfo{address}{Cambridge},
  \bibinfo{year}{2001}).

\bibitem[{\citenamefont{Durbin et~al.}(1998)\citenamefont{Durbin, Eddy, Krogh,
  and Mitchison}}]{Durbin1998}
\bibinfo{author}{\bibfnamefont{R.}~\bibnamefont{Durbin}},
  \bibinfo{author}{\bibfnamefont{S.}~\bibnamefont{Eddy}},
  \bibinfo{author}{\bibfnamefont{A.}~\bibnamefont{Krogh}}, \bibnamefont{and}
  \bibinfo{author}{\bibfnamefont{G.}~\bibnamefont{Mitchison}},
  \emph{\bibinfo{title}{Biological Sequence Analysis}}
  (\bibinfo{publisher}{Cambridge University Press},
  \bibinfo{address}{Cambridge}, \bibinfo{year}{1998}).

\bibitem[{\citenamefont{Cover and Thomas}(1991)}]{Cover1991}
\bibinfo{author}{\bibfnamefont{T.~M.} \bibnamefont{Cover}} \bibnamefont{and}
  \bibinfo{author}{\bibfnamefont{J.~A.} \bibnamefont{Thomas}},
  \emph{\bibinfo{title}{Elements of Information Theory}}
  (\bibinfo{publisher}{Wiley-Interscience}, \bibinfo{address}{New York},
  \bibinfo{year}{1991}).

\bibitem[{\citenamefont{MacKay}(2003)}]{MacKay2003}
\bibinfo{author}{\bibfnamefont{D.~J.~C.} \bibnamefont{MacKay}},
  \emph{\bibinfo{title}{Information Theory, Inference, and Learning
  Algorithms}} (\bibinfo{publisher}{Cambridge University Press},
  \bibinfo{address}{Cambridge}, \bibinfo{year}{2003}).

\bibitem[{\citenamefont{Samengo}(2002)}]{Samengo2002}
\bibinfo{author}{\bibfnamefont{I.}~\bibnamefont{Samengo}},
  \bibinfo{journal}{Phys. Rev. E} \textbf{\bibinfo{volume}{65}},
  \bibinfo{pages}{46124} (\bibinfo{year}{2002}).

\bibitem[{\citenamefont{Young and Crutchfield}(1994)}]{Young1994}
\bibinfo{author}{\bibfnamefont{K.}~\bibnamefont{Young}} \bibnamefont{and}
  \bibinfo{author}{\bibfnamefont{J.~P.} \bibnamefont{Crutchfield}},
  \bibinfo{journal}{Chaos, Solitons, and Fractals}
  \textbf{\bibinfo{volume}{4}}, \bibinfo{pages}{5 } (\bibinfo{year}{1994}).

\bibitem[{\citenamefont{Abramowitz and Stegun}(1965)}]{Abramowitz1965}
\bibinfo{author}{\bibfnamefont{M.}~\bibnamefont{Abramowitz}} \bibnamefont{and}
  \bibinfo{author}{\bibfnamefont{I.~A.} \bibnamefont{Stegun}},
  \emph{\bibinfo{title}{Handbook of Mathematical Functions}}
  (\bibinfo{publisher}{Dover}, \bibinfo{address}{New York},
  \bibinfo{year}{1965}).

\bibitem[{\citenamefont{Crutchfield}(1994)}]{Crutchfield1994}
\bibinfo{author}{\bibfnamefont{J.~P.} \bibnamefont{Crutchfield}},
  \bibinfo{journal}{Physica D} \textbf{\bibinfo{volume}{75}},
  \bibinfo{pages}{11} (\bibinfo{year}{1994}).

\bibitem[{\citenamefont{Upper}(1997)}]{Upper1997}
\bibinfo{author}{\bibfnamefont{D.~R.} \bibnamefont{Upper}}, Ph.D. thesis,
  \bibinfo{school}{University of California}, \bibinfo{address}{Berkeley}
  (\bibinfo{year}{1997}), \bibinfo{note}{{P}ublished by University Microfilms
  Intl, Ann Arbor, Michigan}.

\bibitem[{\citenamefont{Blackwell and Koopmans}(1957)}]{Blackwell1957}
\bibinfo{author}{\bibfnamefont{D.}~\bibnamefont{Blackwell}} \bibnamefont{and}
  \bibinfo{author}{\bibfnamefont{L.}~\bibnamefont{Koopmans}},
  \bibinfo{journal}{Ann. Math. Stat.} \textbf{\bibinfo{volume}{28}},
  \bibinfo{pages}{1011} (\bibinfo{year}{1957}).

\bibitem[{\citenamefont{Wilks}(1962)}]{Wilks1962}
\bibinfo{author}{\bibfnamefont{S.~S.} \bibnamefont{Wilks}},
  \emph{\bibinfo{title}{Mathematical Statistics}} (\bibinfo{publisher}{John
  Wiley \& Sons, Inc.}, \bibinfo{address}{New York}, \bibinfo{year}{1962}).

\bibitem[{\citenamefont{Majumder and
  Bhattacharjee}(1973{\natexlab{a}})}]{Majumder1973}
\bibinfo{author}{\bibfnamefont{K.}~\bibnamefont{Majumder}} \bibnamefont{and}
  \bibinfo{author}{\bibfnamefont{G.}~\bibnamefont{Bhattacharjee}},
  \bibinfo{journal}{Appl. Stat.} \textbf{\bibinfo{volume}{22}},
  \bibinfo{pages}{411} (\bibinfo{year}{1973}{\natexlab{a}}).

\bibitem[{\citenamefont{Majumder and
  Bhattacharjee}(1973{\natexlab{b}})}]{Majumder1973a}
\bibinfo{author}{\bibfnamefont{K.}~\bibnamefont{Majumder}} \bibnamefont{and}
  \bibinfo{author}{\bibfnamefont{G.}~\bibnamefont{Bhattacharjee}},
  \bibinfo{journal}{Appl. Stat.} \textbf{\bibinfo{volume}{22}},
  \bibinfo{pages}{409} (\bibinfo{year}{1973}{\natexlab{b}}).

\bibitem[{\citenamefont{Cran et~al.}(1977)\citenamefont{Cran, Martin, and
  Thomas}}]{Cran1977}
\bibinfo{author}{\bibfnamefont{G.}~\bibnamefont{Cran}},
  \bibinfo{author}{\bibfnamefont{K.}~\bibnamefont{Martin}}, \bibnamefont{and}
  \bibinfo{author}{\bibfnamefont{G.}~\bibnamefont{Thomas}},
  \bibinfo{journal}{Appl. Stat.} \textbf{\bibinfo{volume}{26}},
  \bibinfo{pages}{111} (\bibinfo{year}{1977}).

\bibitem[{\citenamefont{Berry et~al.}(1990)\citenamefont{Berry, {P.W. Mielke,
  Jr.}, and Cran}}]{Berry1990}
\bibinfo{author}{\bibfnamefont{K.}~\bibnamefont{Berry}},
  \bibinfo{author}{\bibnamefont{{P.W. Mielke, Jr.}}}, \bibnamefont{and}
  \bibinfo{author}{\bibfnamefont{G.}~\bibnamefont{Cran}},
  \bibinfo{journal}{Appl. Stat.} \textbf{\bibinfo{volume}{39}},
  \bibinfo{pages}{309} (\bibinfo{year}{1990}).

\end{thebibliography}
\end{document}